\documentclass[12pt]{amsart}

\usepackage{eucal,amsfonts,amsmath}
\usepackage{amsmath}
\usepackage{amssymb}
\usepackage{latexsym}
\usepackage[dvips]{graphicx}
\usepackage{amsfonts}
\usepackage{amssymb}
\usepackage{amsthm}
\usepackage{latexsym}

\newtheorem{thm}{Theorem}[section]
\newtheorem{cor}[thm]{Corollary}
\newtheorem{defn}[thm]{Definition}
\newtheorem{example}[thm]{Example}

\newtheorem{remark}[thm]{Remark}

\numberwithin{equation}{section}

\newcommand{\LL}{{\mathcal L}}
\newcommand{\U}{{\mathcal U}}
\newcommand{\A}{{\mathcal A}}
\newcommand{\E}{{\mathcal E}}

\newcommand{\Z}{\mathbb{Z}}

\newcommand{\C}{\mathbb{C}}
\newcommand{\K}{\mathbb{K}}
\newcommand{\PP}{\mathbb{P}}

\newcommand{\T}{\mathbb{T}}

\begin{document}

\title[Fundamental groups of curve complements]
{Obstructions on fundamental groups\\ of plane curve complements}
\author[Constance Leidy]{Constance Leidy}
\address{C. Leidy: Department of Mathematics,
          University of Pennsylvania,
          209 S 33rd St., Philadelphia, PA, 19104-6395,
                 USA.}
\email{cleidy@math.upenn.edu}

\author[Laurentiu Maxim ]{Laurentiu Maxim}
\address{L. Maxim : Department of Mathematics,
          University of Illinois at Chicago,
          851 S Morgan Street, Chicago, IL 60607,
                 USA.}
\email{lmaxim@math.uic.edu}

\date{\today}

\maketitle

\begin{abstract}We survey various Alexander-type invariants of plane curve
complements, with an emphasis on obstructions on the type of groups
that can arise as fundamental groups of complements to complex plane
curves. Also included are some new computations of higher-order
degrees of curves, which are invariants defined in a previous paper
of the authors.
\end{abstract}

\section{ Introduction }

This paper is an attempt to give partial answers to the following
question posed by Serre: \emph{what restrictions are imposed on a
group by the fact that it can appear as fundamental group of a
smooth algebraic variety?} There are characteristic zero and
respectively finite characteristic aspects of this problem, but we
will restrict ourselves to the zero characteristic case. More
precisely, our ground field will be $\C$. In what follows we only
treat the very special case of open varieties which are complements
to hypersurfaces in $\C^n$ (note that complements to closed
varieties of complex codimension at least two are simply-connected).

By a Zariski theorem of Lefschetz type (see \cite{Di}, Thm. 1.6.5),
for a generic plane $E$ relative to a given hypersurface $V \subset
\C^n$, the natural map
$$\pi_1(E-E \cap V) \to \pi_1 (\C^n -V)$$ is an isomorphism.
Therefore, possible fundamental groups of complements to
hypersurfaces in $\C^n$ are precisely the fundamental groups of
plane affine curve complements.  Thus, it suffices to restrict
ourselves to the case of complements to curves in $\C^2$.

In view of the above, we can ask now the following refinement of
Serre's question: \emph{ what groups can be realized as fundamental
groups of plane curve complements? what obstructions are there?}

In the next sections we will discuss invariants of the fundamental
group of an affine plane curve complement that are obtained by
studying certain covering spaces of the complement: the Alexander
polynomial is an invariant of the total linking number infinite
cyclic cover, characteristic varieties (in particular, the support)
are derived from studying the universal abelian cover, and the
higher-order degrees are numerical invariants obtained by studying
certain solvable covers associated to terms of the rational derived
series of the group. We will see that these invariants obstruct many
knot groups from being realized as fundamental groups of plane curve
complements.

In the last section we include some examples of explicit
calculations of the higher-order degrees associated to some curve
complements. We will also find some examples of groups that cannot
be realized as the fundamental group of a curve complement (in
general position at infinity) because the higher-order degrees
obstruct this.

\section{ Plane curve complements }

Throughout this paper, we consider the following setting: Let $G$ be
a group, and assume that there is a reduced curve $C=\{f(x,y)=0\}$
in $\C^2$ of degree $d$, with $s$ irreducible components, such that
$G=\pi_1(\C^2 \setminus C)$. For simplicity, we assume that $C$ is
in general position at infinity, that is, its projective completion
is transverse to the line at infinity, though many results remain
valid without this restriction on the behavior at infinity.

We will perform the dual task of studying topological properties of
the curve by studying the fundamental group of its complement, while
at the same time deriving obstructions on a group imposed by the
fact that it is the fundamental group of an affine plane curve
complement. For more comprehensive surveys on the topology of plane
curves and a list of open problems, the interested reader may also
consult the papers \cite{Li15, Li16, O}.

First note that $H_1(G)=H_1(\C^2 \setminus C)=G/G'=\Z^s$, generated
by meridians about the smooth parts of irreducible components of
$C$.

Although in geometric problems fundamental groups of complements to
projective curves play a central role, by switching to the affine
setting (i.e., by also removing a generic line) no essential
information is lost. Indeed, if $\bar C \subset \C\PP^2$ is the
projective completion of $C$, the two groups are related by the
central extension \begin{equation}\label{les} 0 \to \Z \to
\pi_1(\C\PP^2-(\bar{C} \cup H)) \to \pi_1(\C\PP^2 - \bar {C}) \to
0.\end{equation} Moreover, by \cite{O}, Lemma 2, the commutator
subgroups of the affine and respectively projective complements
coincide: \begin{equation} G'=\pi_1(\C\PP^2 -
\bar{C})'.\end{equation}

\subsection{ The linking number infinite cyclic cover of the complement }\label{ic}
We begin with a brief survey of results on the Alexander polynomial
of the curve $C$.

Let $lk:G=\pi_1(\C^2-C) \to \Z$ be the total linking number
epimorphism, i.e. $\alpha \mapsto lk\#(\alpha,C)$. Note that $lk$
factors through $H_1(G)$, sending the basis vectors of $\Z^s$ to
$1$. Let $\U^c$ be the covering of $\U$ corresponding to $Ker(lk)$.
$\U^c$ will be called the \emph{total linking number infinite cyclic
cover of the complement}.

The group of deck transformations of $\U^c$ is $\Z$, and it acts on
$H_1(\U^c;\C)$ by a generating transformation, thus making
$H_1(\U^c;\C)$ into a module over $\C[\Z]=\C[t,t^{-1}]$. This module
is called \emph{the infinite cyclic Alexander module of the curve
complement}. As $\C[t,t^{-1}]$ is a principal ideal domain,
$H_1(\U^c;\C)$ decomposes as
$$H_1(\U^c;\C) \cong \C[t,t^{-1}]^m \oplus \left( \oplus_i \C[t,t^{-1}]/\lambda_i(t)
\right),$$ for some $m \in \Z$ and polynomials $\lambda_i(t)$
defined up to a unit of $\C[t,t^{-1}]$. In fact, the following
result holds:

\begin{thm} (Zariski-Libgober \cite{Li0}) $H_1(\U^c;\C)$ is a torsion
$\C[t,t^{-1}]$-module.
\end{thm}

Therefore, it does make sense to associate to $C$ a polynomial,
namely the order of $H_1(\U^c;\C)$ (cf. \cite{Mi2}). This is a
global invariant of $C$ (or of $G$) defined as follows:
\begin{defn}\rm $\Delta_C(t)=\prod_i \lambda_i(t)$ is
called \emph{the Alexander polynomial of $C$ (or $G$)}.
\end{defn}

It is easy to see that the exponent of $(t-1)$ in $\Delta_C(t)$ is
$s-1$, where $s$ is the number of irreducible components of $C$
(e.g., see \cite{O}). In particular, if the curve $C$ is
irreducible, the Alexander polynomial $\Delta_C(t)$ can be
normalized so that $\Delta_C(1)=1$ (cf. \cite{Li0}).

\subsubsection{ Libgober's divisibility theorem for Alexander polynomials }
In \cite{Li0, Li2, Li22}, Libgober gives an algebraic-geometrical
meaning of the Alexander polynomial of $C$ as follows.

With each singular point $x \in C$ there is an associated
\emph{local Alexander polynomial}, $\Delta_x(t)$, defined as the
characteristic polynomial of the monodromy of local Milnor fibration
at $x$ (cf. \cite{Mi}). Then:

\begin{thm}\label{div} (Libgober \cite{Li0, Li2}) Up to a power of $(t-1)$, the Alexander polynomial
$\Delta_C(t)$ of a plane curve in general position at infinity
divides the product $\prod_{x \in Sing(C)} \Delta_x(t)$ of the local
Alexander polynomials at the singular points of $C$. Therefore the
local type of singularities has an effect on the topology of $C$.
\end{thm}

Zariski also showed that \emph{the position of singularities} has an
influence on the topology of $C$. Moreover, as Libgober observed,
the Alexander polynomial is sensitive to the position of
singularities (\cite{Li0}). The classical example of Zariski's
sextics with six cusps will be discussed in section \ref{ex}.

Theorem \ref{div} remains true without any assumption on the
behavior of $C$ at infinity, but one has to take into account the
contribution of singularities at infinity. As a corollary of this
fact, we have that

\begin{cor} $\Delta_C(t)$ is cyclotomic. Moreover,
for a curve $C$ in general position at infinity, the zeros of
$\Delta_C(t)$ are roots of unity of order $d=deg(C)$.
\end{cor}

It follows that many knot groups, e.g. that of figure eight knot
(whose Alexander polynomial is $t^2-3t+1$), cannot be of the form
$\pi_1(\C^2-C)$. However, the class of possible fundamental groups
of plane curve complements includes braid groups, or groups of torus
knots of type $(p,q)$ (see \cite{LM} \S 5, and the references
therein).

\begin{remark}\label{M}\rm The above divisibility result has been
generalized to higher dimensions by Libgober (\cite{Li1, Li}), who
considered complements to affine hypersurfaces with only isolated
singularities, and also by Maxim (\cite{M}), who in his thesis
treated the case of hypersurfaces with non-isolated singularities.
From a divisibility result in \cite{M}, it follows that in
Libgober's divisibility result it suffices to consider only the
contribution of local Alexander polynomials at singular points
contained in some fixed irreducible component of the hypersurface.
In particular, this shows that the Alexander polynomial does not
provide enough information about the topology of reducible curves
(hypersurfaces). For example, if $C$ is a union of two curves that
intersect transversally, then $\Delta_C(t)=(t-1)^{s-1}$  (see
\cite{O}). To overcome this problem, we study higher coverings of
the complement.
\end{remark}

\subsection{ The universal abelian cover of the complement }\label{ua}
In this section, following \cite{Li4} we define invariants
associated to the universal abelian cover of the complement.

Let $\U^{ab}$ be the universal abelian cover of $\U$, i.e., the
covering associated to the subgroup $G'$. Under the action of the
covering transformation group, the \emph{universal abelian module}
$H_1(\U^{ab};\C)=G'/G'' \otimes \C$ becomes a finitely generated
module over $\C[G/G']=\C[t_1^{\pm 1},..,t_s^{\pm 1}]=:R_s$. Note
that $R_s$ is a Noetherian domain and a UFD.

Now let $M$ be a presentation matrix of $\A:=H_1(\U^{ab};\C)$
corresponding to a sequence $$(R_s)^m \to (R_s)^n \to \A \to 0$$

\begin{defn}\rm The \emph{order ideal} of $\A$, $\E_0(\A)$, is the
ideal in $R_s$ generated by the $n \times n$-minor determinants of
$M$, with the convention $\E_0(\A)=0$ if $n > m$. The \emph{support}
of $\A$, $Supp(\A)$, is the reduced sub-scheme of the s-dimensional
torus $\T^s=Spec(R_s)$ defined by the order ideal. Equivalently, a
prime ideal $p$ is in $Supp(\A)$ if and only if $\A_p \neq 0$ (that
is, if and only if $p \supset Ann(\A)$).
\end{defn}

\noindent Similarly, the \emph{$i$-th (algebraic) characteristic
variety} is defined by the $i$-th elementary ideal of $\A$. Away
from the trivial character, characteristic varieties of $\A$
coincide with jumping loci of homology of rank-one local systems on
the complement (cf. \cite{Li3}), defined as
$$V^t_i(G)=\{\lambda \in {\C^{*}}^s | dim_{\C}H_1(G, \LL_{\lambda}) \geq i \}, \ 1 \leq i \leq s,$$
where $\LL_{\lambda}$ is the rank-one local system associated to the
character $\lambda$. In \cite{DM}, these jumping loci are called
\emph{topological characteristic varieties}. By a result of Arapura
(\cite{A}), each $V^t_i(G)$ is a union of subtori of the character
torus, possibly translated by unitary characters. This fact imposes
strong obstructions on the group $G$. Characteristic varieties, both
algebraic and topological, give very precise information about the
homology of (finite) abelian covers of $\U$ (e.g., see \cite{Li4}).

\begin{example}\rm
(1) \ If $C$ is irreducible, then
$Supp(\A)=\{\Delta_C(t)=0\}$.\newline (2) \ If $L$ is a link in
$S^3$ and $G=\pi_1(S^3-L)$ then $Supp(\A)$ is the zero-set of the
multivariable Alexander polynomial of the link.
\end{example}

\begin{remark}\rm A multivariable Alexander polynomial of $C$
could be defined as the greatest common divisor of all elements of
the order ideal $\E_0(\A)$. However, if $\text{codim}_{\T^s}
Supp(\A)
>1$, then this polynomial is trivial, so it doesn't contain any interesting information about the topology of $C$.
\end{remark}

The support of the universal abelian module is restricted by the
following result (\cite{Li4}, \cite{DM}):
\begin{thm} (Libgober) If $C$ is a curve in general
position at infinity, then
$$Supp(A) \subset \{ (\lambda_1,..,\lambda_s)\in \T^s \ \ | \ \
\prod_{i=1}^s \lambda_i ^{d_i}=1 \}$$ where $d_i$ is the degree of
the $i$-th irreducible component of $C$.
\end{thm}

In \cite{DM}, there is a similar characterization of supports of
universal abelian invariants associated to complements of
hypersurfaces in $\C^{n+1}$, with any type of singularities. The
supports are also shown to depend on the local type of
singularities.

\subsection{ Higher-order coverings of the complement }
In this section we study covers of the curve complement that are
associated to terms in the rational derived series of the
fundamental group. The invariants arising in this way were
originally used in the study of knots and respectively
$3$-manifolds, e.g. to show that certain groups cannot be realized
as the fundamental group of the complement of a knot, or as the
fundamental group of a 3-manifold. Some very useful background
material is presented in \cite{C,H}.

Let $G_r ^{(0)}=G$. For $n \geq 1$, we define the $n^{th}$ term of
the \emph{rational derived series of $G$} inductively by:
$$G_r ^{(n)}=\{g \in G_r ^{(n-1)} | g^k \in
[G_r ^{(n-1)},G_r ^{(n-1)}], \ \text{for some} \ k \in \Z -\{0\}
\}.$$ It is easy to see that $G_r ^{(i)} \triangleleft G_r ^{(j)}
\triangleleft G$, if $i \geq j \geq 0$, so we can consider quotient
groups. Set $\Gamma_n:=G/G_r ^{(n+1)}$. We use rational derived
series as opposed to the usual derived series in order to avoid
zero-divisors in the group ring $\Z\Gamma_n$.

The successive quotients of the rational derived series are
torsion-free abelian groups. Indeed (cf. \cite{H}),
$$G_r^{(n)}/G_r^{(n+1)} \cong \left(G_r^{(n)}/[G_r^{(n)},
G_r^{(n)}] \right)/\{\Z-\text{torsion}\}.$$ Therefore, if
$G=\pi_1(\C^2-C)$, then $G'=G_r'$ (this follows from the trivial
fact that $G'$ is a subgroup of $G'_r$, together with $G/G' \cong
\Z^s$).

By construction, it follows that $\Gamma_n$ is a
poly-torsion-free-abelian group, in short PTFA (\cite{H}), i.e., it
admits a normal series of subgroups such that each of the successive
quotients of the series is torsion-free abelian. Then $\Z\Gamma_n$
is a right and left Ore domain, so it embeds in its classical right
ring of quotients $\mathcal{K}_n$, a skew-field.

\begin{defn}\rm The \emph{$n$-th order Alexander modules of $C$}
are
$$\mathcal{A}^{\Z}_n(C)=H_1(\U;\Z\Gamma_n)=H_1(\U_{\Gamma_n};\Z)$$
where $\U_{\Gamma_n}$ is the covering of $\U$ corresponding to the
subgroup $G_r^{(n+1)}$.  That is,
$\mathcal{A}^{\Z}_n(C)=G_r^{(n+1)}/[G_r^{(n+1)},G_r^{(n+1)}]$ as a
right $\Z \Gamma_n$-module.\newline The \emph{$n^{th}$ order rank}
of (the complement of) $C$ is:
$$r_n(C)=\text{rk}_{\mathcal{K}_n}H_1(\U;\mathcal{K}_n)$$
\end{defn}

\begin{remark}\rm
Note that
$\mathcal{A}^{\Z}_0(C)=G_r^{(1)}/[G_r^{(1)},G_r^{(1)}]=G'/G''$ is
just the universal abelian invariant of the complement.
\end{remark}

\begin{remark}\rm If $C$ is an irreducible curve (or
$\beta_1(G)=1$), it follows directly from a result in \cite{C} that
$\mathcal{A}^{\Z}_n(C)$ is a torsion $\Z\Gamma_n$-module. In
\cite{LM}, the authors showed that this is also true for the
reducible case (at least for curves in general position at
infinity). (See Theorem \ref{f}.)
\end{remark}

\begin{example}\label{sm-nod}\rm (1) \ If $C$ is
non-singular and in general position at infinity, then $G=\Z$.
\newline (2) \ If $C$ has only nodal singular points (locally defined
by $x^2-y^2=0$), then $G$ is abelian. \newline In both cases above
it follows that $\mathcal{A}^{\Z}_0(C)=0$, and therefore
$\mathcal{A}^{\Z}_n(C)=0$ for all $n$ (cf. \cite{LM}, Remark
3.4).\end{example}

We associate to any curve $C$ (or equivalently, to its group $G$) a
sequence of non-negative integers $\delta_n(C)$ as follows (it is
more convenient to work over a PID, so we look for a ``convenient"
one): Let $\psi \in H^1(G;\Z)$ be the primitive class representing
the linking number homomorphism $G \overset{\psi}{\to} \Z$, $\alpha
\mapsto \text{lk}(\alpha,C)$. Since $G'$ is in the kernel of $\psi$,
we have a well-defined induced epimorphism $\bar{\psi} : \Gamma_n
\to \Z$. Let $\bar{\Gamma}_n =Ker \bar{\psi}$. Then $\bar{\Gamma}_n$
is a PTFA group, so $\Z\bar{\Gamma}_n$ has a right ring of quotients
$\K_n=(\Z\bar{\Gamma}_n)S_n ^{-1}$, for $S_n=\Z\bar{\Gamma}_n - 0$.
Set $R_n=(\Z\Gamma_n)S_n^{-1}$. Then $R_n$ is a flat left
$\Z\Gamma_n$-module.

Crucially, $R_n$ is a PID, isomorphic to the ring of skew-Laurent
polynomials $\K_n[t^{\pm 1}]$. Indeed, by choosing a $t \in
\Gamma_n$ such that $\bar {\psi} (t)=1$, we get a splitting $\phi$
of $\bar{\psi}$, and the embedding $\Z\bar{\Gamma}_n \subset \K_n$
extends to an isomorphism $R_n \cong \K_n[t^{\pm 1}]$. However this
isomorphism depends in general on the choice of splitting of $\bar
\psi$ !

\begin{defn}\rm
(1) \ The \emph{$n^{th}$-order localized Alexander module of the
curve $C$} is defined to be $\mathcal{A}_n(C)=H_1(\U;R_n)$, viewed
as a right $R_n$-module.  If we choose a splitting $\phi$ to
identify $R_n$ with $\K_n[t^{\pm 1}]$, we define
$\mathcal{A}^{\phi}_n(C)=H_1(\U;\K_n[t^{\pm 1}])$. \newline (2) \
The \emph{$n^{th}$-order degree of $C$} is defined to be:
$$\delta_n(C)=\text{rk}_{\K_n} \mathcal{A}_n(C)=\text{rk}_{\K_n} \mathcal{A}^{\phi}_n(C).$$
\end{defn}

\begin{remark}\rm Note that $\delta_n(C) < \infty$ if and only if
$\text{rk}_{\mathcal{K}_n}H_1(\U;\mathcal{K}_n)=0$, i.e.
$\mathcal{A}_n(C)$ is a torsion module.
\end{remark}

The degrees $\delta_n(C)$ are integral invariants of the fundamental
group $G$ of the complement. Indeed, by \cite{Har} \S 1, we have:
$$\delta_n(C)=rk_{\K_n} \left(
{G^{(n+1)}_r}/[G^{(n+1)}_r,G^{(n+1)}_r] \otimes_{\Z\bar{\Gamma}_n}
\K_n \right).$$ Since the isomorphism between $R_n$ and $\K_n[t^{\pm
1}]$ depends on the choice of splitting, we \emph{cannot} define in
a meaningful way a ``higher-order Alexander polynomial", as we did
in the infinite cyclic case. However, for any \emph{choice} of
splitting, the degree of the associated higher-order Alexander
polynomial is the same. Therefore although a higher-order Alexander
polynomial is not well-defined in general, the \emph{degree} of the
polynomial associated to a choice of a splitting yields a
well-defined invariant of $G$. This is exactly the higher-order
degree $\delta_n$ defined above.

The higher-order degrees of $C$ may be computed by means of Fox free
calculus by using a presentation of $\pi_1(\C^2-C)$. The latter can
be obtained by means of Moishezon's braid monodromy \cite{Mo}. In
general, these steps are difficult to achieve. However in section
\ref{ex}, some examples are explicitly computed.

The obstructions on $G$ obtained from analyzing the higher-order
degrees of a plane curve complement are contained in the following
theorem of \cite{LM}:

\begin{thm}\label{f} (Leidy-Maxim \cite{LM})\newline If $G=\pi_1(\C^2-C)$ for some plane curve $C$ in
general position at infinity, then the higher-order degrees
$\delta_n(C)$ are finite. More precisely:
\begin{enumerate}
\item there exists a uniform upper bound in terms of the degree of
$C$: $\delta_n(C) \leq d(d-2)$, for all $n$.
\item for each $n$, there is an upper bound in terms of local
invariants at singular points of $C$
$$\delta_n(C) \leq \Sigma_{k=1}^l \left(\mu(C,c_k) + 2
n_k\right)+ 2g + d - l$$ where $c_k$, $1 \leq k \leq l$, are the
singularities of $C$, $n_k$ is the number of branches through the
singularity $c_k$, $\mu(C,c_k)$ is the Milnor number of the
singularity germ $(C,c_k)$, and $g$ is the genus of the normalized
curve.
\end{enumerate}
\end{thm}

We have the following important corollary that provides an
obstruction to a group being the fundamental group of the complement
of a curve in general position at infinity. This can be combined
with the central extension (\ref{les}) in order to obtain
obstructions on the fundamental groups of projective plane curve
complements. (In the last section of this paper we will use this
corollary to find such examples.)

\begin{cor} If $C$ is a plane curve in general
position at infinity, then $\mathcal{A}^{\Z}_n(C)$ is a torsion
$\Z\Gamma_n$-module.
\end{cor}

\section{ Examples }\label{ex}

In this section, we will present some explicit calculations of the
higher-order degrees of various curve complements. Although
computing higher-order degrees can be difficult, we hope that these
examples will aide the reader in understanding how a general
calculation can be carried out.

Before presenting the calculations, we recall some results from
\cite{LM}.

\begin{itemize}

\item If $C$ is either non-singular or has only nodal singular points
(and is in general position at infinity), it follows from Example
\ref{sm-nod} that $\delta_n(C)=0$ for all $n\geq0$.

\item If $C$ is defined by a weighted
homogeneous polynomial $f(x,y)=0$, then either:
\begin{itemize}
\item if either $n>0$ or $\beta_1(\U)>1$, then
$\delta_n(C)=\mu(C,0)-1$.
\item if $\beta_1(\U)=1$, then $\delta_0(C)=\mu(C,0)$, where
$\mu(C,0)$ is the Milnor number of the singularity germ at the
origin.
\end{itemize}

\item If $C$ is an irreducible affine curve, then
$\delta_0(C)=\text{deg} \Delta_C(t)$, where $\Delta_C(t)$ denotes
the Alexander polynomial of the curve complement. If, moreover, the
Alexander polynomial is trivial then all higher-order degrees
vanish.
\end{itemize}

In the next three examples, we consider an irreducible curve $\bar C
\subset \C\PP^2$ and a generic line (at infinity) $H$, then set
$C=\bar C -H$.

\begin{example}\rm Let $\bar C \subset \C\PP^2$ be a degree
$d$ curve having only nodes and cusps as its only singularities. If
$d \not\equiv 0 \ (\text{mod} \ 6)$, then all higher-order degrees
of $C$ vanish. (this follows from the divisibility results on
$\Delta_C(t)$, which imply that $\Delta_C(t)=1$).
\end{example}

\begin{example}\rm If $\bar C$ is Zariski's three-cuspidal
quartic, then $G=\pi_1(\C^2-C)=\langle a, b \ | \ aba=bab, a^2=b^2
\rangle.$ Thus $G' \cong \Z/3\Z$. So $\delta_n(C)=0$, for all $n$.
For all other quartics, the corresponding group of the affine
complement is abelian, so the higher-order degrees vanish again.
\end{example}

\begin{example}\rm \emph{Zariski's sextics with $6$ cusps} \newline
Let $\bar C \subset \C\PP^2$ be a curve of degree $6$ with
$6$ cusps.
\begin{itemize} \item If the $6$ cusps are on a conic, then
$\pi_1(\C^2 - C)=\pi_1(\C\PP^2- \bar C \cup H)$ is isomorphic to the
fundamental group of the trefoil knot, and has Alexander polynomial
$t^2-t+1$. Thus, $\delta_0(C)=2$, and $\delta_n(C)=1$ for all $n>0$.
\item If the six cusps are not on a conic, then
$\pi_1(\C^2- C)$ is abelian. Therefore, $\delta_n(C)=0$ for all $n
\geq 0$.
\end{itemize}
\end{example}

\begin{remark}\rm From the above example we see that the higher-order degrees of
a curve, at any level $n$, are also sensitive to the position of
singular points. An interesting open problem is to find Zariski
pairs that are distinguished by some $\delta_k$, but not
distinguished by any $\delta_n$ for $n < k$.
\end{remark}

\subsection{Line Arrangements}
Since we are assuming that our curves are in generic position at
infinity, the arrangements that we will consider do not have
parallel lines. If we have an arrangement with two intersecting
lines, the only singularity is a node, and therefore $\delta_n$ is
trivial for all $n$. Similarly, $\delta_n=0$ if we have three lines
arranged so that the singularities are each nodes. Hence the first
interesting case to consider is the arrangement of three lines
intersecting in a triple point. Using the techniques of \cite{CS} we
can find a presentation for the fundamental group of $\U$, the
complement of the three lines in $\C^2$:
$$\pi_1(\U)\cong\langle \sigma_1,\sigma_2,\sigma_3|\sigma_1\sigma_2\sigma_3
=\sigma_2\sigma_3\sigma_1=\sigma_3\sigma_1\sigma_2\rangle.$$ Here
$\sigma_1$, $\sigma_2$, and $\sigma_3$ correspond to the meridians
of the lines. In particular, they each map to a different generator
of of $H_1(\U)\cong\Z^3$ and they are all mapped to the same
generator of $\Z$ under the total linking number homomorphism
$\pi_1(\U) \to H_1(\U) \to \Z$. It is easier to work with a
presentation for $\pi_1(\U)$ where only one generator maps to the
generator of $\Z$ under the total linking number homomorphism. Hence
we choose new generators: $a=\sigma_1$, $b=\sigma_2\sigma_1^{-1}$,
and $c=\sigma_3\sigma_1^{-1}$. With these new generators we have the
following presentation:
$$\pi_1(\U)\cong\langle a,b,c|abac=baca=ca^2b\rangle.$$

Using Fox calculus \cite{F1}, \cite{F2}, we can obtain a
presentation matrix for $H_1(\U,u_0;\Z\pi_1(\U))$, the homology of
the universal cover of $\U$ relative to a basepoint $u_0$ as a left
$\Z\pi_1(\U)$-module. The 1-chains for the universal cover of $\U$
are generated as a $\Z\pi_1(\U)$-module by $\alpha$, $\beta$, and
$\gamma$, where $\alpha$, $\beta$, and $\gamma$ each represent a
single lift of the 1-chains of $\U$ corresponding to $a$, $b$, and
$c$, respectively. Since we are computing the homology relative to a
basepoint, $\alpha$, $\beta$, and $\gamma$ are in fact 1-cycles in
$H_1(\U,u_0;\Z\pi_1(\U))$. It remains to consider the 2-chains of
the universal cover of $\U$. First, the 2-chain of $\U$
corresponding to the relation $abaca^{-1}c^{-1}a^{-1}b^{-1}$ in
$\pi_1(\U)$ lifts to a 2-chain of the universal cover of $\U$ whose
boundary is:
\begin{eqnarray*}
&& \alpha + a*\beta + ab*\alpha + aba*\gamma - abaca^{-1}*\alpha -
abaca^{-1}c^{-1}*\gamma \\ && \hspace{.25in} -
abaca^{-1}c^{-1}a^{-1}*\alpha - abaca^{-1}c^{-1}a^{-1}b^{-1}*\beta
\end{eqnarray*}
Using the relation $abaca^{-1}c^{-1}a^{-1}b^{-1}=1$ in $\pi_1(\U)$,
we can rewrite this boundary as:
\begin{eqnarray*}
\alpha + a*\beta + ab*\alpha + aba*\gamma - bac*\alpha - ba*\gamma
-b*\alpha - \beta \\
= (1 + ab -bac - b)*\alpha + (a - 1)*\beta + (aba - ba)*\gamma
\end{eqnarray*}
Similarly, the 2-chain of $\U$ corresponding to the relation
$bacab^{-1}a^{-2}c^{-1}$ in $\pi_1(\U)$ lifts to a 2-chain of the
universal cover of $\U$ whose boundary is:
\begin{eqnarray*}
&&\beta + b*\alpha + ba*\gamma + bac*\alpha - bacab^{-1}*\beta -
bacab^{-1}a^{-1}*\alpha \\ && \hspace{.25in} -
bacab^{-1}a^{-2}*\alpha - bacab^{-1}a^{-2}c^{-1}*\gamma.
\end{eqnarray*}
Using the relation
$bacab^{-1}a^{-2}c^{-1}=1$ in $\pi_1(\U)$, we can rewrite this
boundary as:
\begin{eqnarray*}
\beta + b*\alpha + ba*\gamma + bac*\alpha - ca^2*\beta - ca*\alpha -
c*\alpha - \gamma \\
= (b + bac - ca - c)*\alpha + (1-ca^2)*\beta + (ba - 1)*\gamma
\end{eqnarray*}
We can collect this information to write a presentation matrix:
\begin{equation*}
H_1(\U,u_0;\Z\pi_1(\U))=\left(
\begin{array}{ccc}
1+ab-bac-b & a-1 & aba-ba \\
b+bac-ca-c & 1-ca^2 & ba-1 \\
\end{array}
\right)
\end{equation*}
Here the columns correspond to the generators, $\alpha$, $\beta$,
and $\gamma$, respectively, and the rows correspond to relations.

If we allow elements of $\pi_1(\U)^{(n+1)}_r$ to be set equal to 1
in $\Z\pi_1(\U)$, we can also consider the above as a presentation
matrix for $H_1(\U,u_0;\Z\Gamma_n)$. Furthermore, since $R_n$ is a
flat $\Z\Gamma_n$-module, we can also consider it to be a
presentation matrix for $H_1(\U,u_0;R_n)$. If we think of the matrix
in this way, any non-zero element in $\Z\bar{\Gamma}_n$ has an
inverse. (Recall that $\bar{\Gamma}_n$ is the kernel of the map
$\bar{\psi}:\Gamma_n \to \Z$, induced by the total linking number
homomorphism.)

If we choose a splitting of $\bar{\psi}$, there is an isomorphism
between $R_n$ and $\K_n[t^{\pm1}]$. For our example, we choose the
splitting that maps $t$ to $a$. To obtain a presentation for
$H_1(\U,u_0;\K_n[t^{\pm1}])$ we must replace each entry in the above
matrix with its image under the isomorphism $R_n \to
\K_n[t^{\pm1}]$. This results in the following presentation matrix
for $H_1(\U,u_0;\K_n[t^{\pm1}])$:
\begin{equation*}
\left(
\begin{array}{ccc}
1+aba^{-1}t-baca^{-1}t-b & t-1 & aba^{-1}t^2-bt \\
b+baca^{-1}t-ct-c & 1-ct^2 & bt-1 \\
\end{array}
\right)
\end{equation*}
Notice that because $\K_n[t^{\pm1}]$ is a \emph{skew} Laurent
polynomial ring, we must be careful when writing elements where $t$
is not originally on the right. For example, $tb=aba^{-1}t$ in
$\K_n[t^{\pm1}]$.

The next step in finding $\delta_n$ is diagonalizing this matrix,
which is possible since $\K_n[t^{\pm1}]$ is a PID. Since $c\neq1$ in
$\pi_1(\U)/\pi_1(\U)^\prime$, it follows that $c \notin
\pi_1(\U)^{(n)}_r$ for all $n \geq 1$. Therefore $c\neq1$ in
$\Gamma_n$ for all $n \geq 0$. Hence $1-c\neq0$ in $\Z\Gamma_n$ and
is therefore invertible in $\K_n[t^{\pm1}]$. This allows us to
multiply the last column in our presentation matrix by the unit
$1-c$. Since our matrix is a presentation of a left module and since
columns correspond to generators, we multiply columns on the right.
The result of multiplying the last column (on the right) by the unit
$1-c$ is the following:
\begin{equation*}
\left(
\begin{smallmatrix}
(aba^{-1}-baca^{-1})t+(1-b) & t-1 & (aba^{-1}-abaca^{-2})t^2+(baca^{-1}-b)t \\
(baca^{-1}-c)t+(b-c) & 1-ct^2 & (b-baca^{-1})t+(c-1) \\
\end{smallmatrix}
\right)
\end{equation*}
Next we add the first column times $1-t$ and the second column times
$1-b$ to the last column. The result is the following:
\begin{equation*}
\left(
\begin{array}{ccc}
(aba^{-1}-baca^{-1})t+(1-b) & t-1 & 0 \\
(baca^{-1}-c)t+(b-c) & 1-ct^2 & 0 \\
\end{array}
\right)
\end{equation*}
This means that we have a free generator, which is expected since we
are computing the homology relative to a basepoint.

Next we multiply the first row by $ct+c$ and add it to the second.
Since our matrix is a presentation of a left module and since rows
correspond to relations, we multiply rows on the left. The result of
multiplying the first row (on the left) by $ct+c$ and adding it to
the second is:
\begin{equation*}
\left(
\begin{array}{ccc}
(aba^{-1}-baca^{-1})t+(1-b) & t-1 & 0 \\
(1-c)(baca^{-1}t^2+baca^{-1}t+b) & 1-c & 0 \\
\end{array}
\right)
\end{equation*}
Now we multiply the second row (on the left) by the unit
$(1-c)^{-1}$. Then we multiply the second row by $1-t$ and add it to
the first. This results in the following matrix:
\begin{equation*}
\left(
\begin{array}{ccc}
1-baca^{-1}t^3 & 0 & 0 \\
(baca^{-1}t^2+baca^{-1}t+b) & 1 & 0 \\
\end{array}
\right)
\end{equation*}
Notice that we can now eliminate the second column and row. Hence we
have shown that $H_1(\U,u_0;\K_n[t^{\pm1}])\cong\K_n[t^{\pm1}]
\oplus \K_n[t^{\pm1}]/\langle 1-baca^{-1}t^3\rangle$. To find
$H_1(\U;\K_n[t^{\pm1}])$, we consider the long exact sequence of a
pair:
$$0 \to H_1(\U;\K_n[t^{\pm1}]) \to H_1(\U,u_0;\K_n[t^{\pm1}]) \to H_0(u_0;\K_n[t^{\pm1}]).$$
Since $H_1(\U;\K_n[t^{\pm1}])$ is a torsion module and
$H_0(u_0;\K_n[t^{\pm1}])$ is a free module, we conclude that
$H_1(\U,u_0;\K_n[t^{\pm1}])\cong\K_n[t^{\pm1}]/\langle
1-baca^{-1}t^3\rangle$. Therefore, for the arrangement of three
lines intersecting in a triple point, $\delta_n=3$ for all $n\geq0$.

If we add an additional line to this arrangement that intersects
previous three lines in nodes (as in the wiring diagram below),
$\delta_n=0$ for all $n\geq0$.
\begin{figure}[h]
\centering
\includegraphics[scale=1.5]{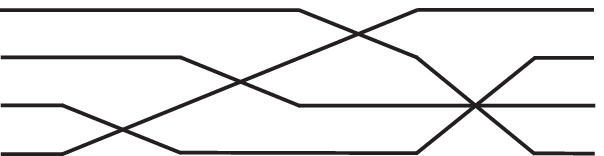}
\end{figure}
In fact, for any line arrangement that contains a line whose only
intersections are nodes, $\delta_n=0$ for all $n\geq0$.

If instead we add an additional line to this arrangement so that all
lines intersect in a single point, $\delta_n=8$ for all $n\geq0$.
The arrangement of five lines intersecting in a single point has
$\delta_n=15$ for all $n\geq0$. Each of these calculations can be
done in the same fashion as the one above. We conjecture that for m
lines intersecting in a single point, $\delta_n=m(m-2)$ for all
$n\geq0$.

\subsection{Artin groups of spherical-type}
Deligne \cite{D} showed that each Artin group of spherical-type
appears as the fundamental group of the complement of a complex
hyperplane arrangement. Mulholland and Rolfsen \cite{MR} showed that
the commutator subgroups of the following Artin groups are perfect
(i.e. $G'=G''$): $A_n$, $n\geq4$; $B_n$, $n\geq5$; $D_n$, $n\geq5$;
$E_n$, $n=6,7,8$; $H_n$, $n=3,4$. It follows that all higher-order
degrees are trivial for curves whose complements have the above
fundamental groups. We will explicitly compute the higher-order
degrees of the Artin group of type $A_3$.

The Artin group of type $A_3$ is the braid group on four strands,
$\mathfrak{B}_4$. A standard presentation for the braid group on
four strands is:
$$\mathfrak{B}_4=\langle\sigma_1,\sigma_2,\sigma_3|\sigma_1\sigma_3=\sigma_3\sigma_1,
\sigma_1\sigma_2\sigma_1=\sigma_2\sigma_1\sigma_2,
\sigma_2\sigma_3\sigma_2=\sigma_3\sigma_2\sigma_3\rangle.$$ We give
a new presentation by choosing new generators: $x=\sigma_1$,
$y=\sigma_2\sigma_1^{-1}$, and $z=\sigma_3\sigma_1^{-1}$.
$$B_4=\langle x,y,z|xz=zx,xyx=yx^2y,yxzxy=zxyxz\rangle$$
Notice that the abelianization of $\mathfrak{B}_4$ is $\Z$, and that
under the abelianization map, $x$ maps to a generator of $\Z$, while
$y$ and $z$ are mapped to $0$.

If $\U$ is a curve complement with $\pi_1(\U)\cong\mathfrak{B}_4$,
we can use Fox calculus to obtain the following presentation for
$H_1(\U,u_0;\Z\pi_1(\U))$, as a left $\Z\pi_1(\U)$-module:
\begin{equation*}
\left(
\begin{array}{ccc}
1-z & 0 & x-1 \\
1+xy-yx-y & x-yx^2-1 & 0 \\
y+yxz-zxy-z & 1+yxzx-zx & yx-zxyx-1 \\
\end{array}
\right)
\end{equation*}
Here the columns correspond to generators and the rows correspond to
relations.

To obtain a presentation for $H_1(\U,u_0;\K_n[t^{\pm1}])$, we choose
the splitting that maps $t$ to $x$. Then we have the following
presentation for $H_1(\U,u_0;\K_n[t^{\pm1}])$:
\begin{equation*}
\left(
\begin{array}{ccc}
1-z & 0 & t-1 \\
1+xyx^{-1}t-yt-y & t-yt^2-1 & 0 \\
y+yzt-zxyx^{-1}t-z & 1+yzt^2-zt & yt-zxyx^{-1}t^2-1 \\
\end{array}
\right)
\end{equation*}
We remind the reader that $x$ and $z$ commute in $\mathfrak{B}_4$
and therefore in $\K_n[t^{\pm1}]$, we have $tz=zt$.

It follows from Thm 3.6 of \cite{MR} that
$\mathfrak{B}_4^\prime/\mathfrak{B}_4^{\prime\prime}\cong\Z^2$,
generated by $y$ and $xyx^{-1}$. In particular, $y \notin
(\mathfrak{B}_4)^{(n)}_r$ for $n\geq2$. Therefore, $1-y \neq 0$ in
$\Z\Gamma_n$ for $n\geq1$. Hence $1-y$ is invertible in $\K_n$ for
$n\geq1$. We first consider the case when $n=0$ and then continue
the calculation for $n\geq1$.

If $n=0$, then we set $y=z=1$ in the above matrix to obtain:
\begin{equation*}
\left(
\begin{array}{ccc}
0 & 0 & t-1 \\
0 & -t^2+t-1 & 0 \\
0 & t^2-t+1 & -t^2+t-1 \\
\end{array}
\right)
\end{equation*}
After adding the second row along with $t$ times the first row to
the last row, we are able to eliminate the last row and column.
Therefore,
$H_1(\U,u_0;\K_0[t^{\pm1}])\cong\K_0[t^{\pm1}]\oplus\K_0[t^{\pm1}]/\langle
t^2-t+1\rangle$. Hence $\delta_0=2$. (Note that this is simply the
computation of the degree of the classical Alexander polynomial.)

We now assume that $n\geq1$, and therefore can use the fact that
$1-y$ is a unit in $\K_n[t^{\pm1}]$. We begin the process of
diagonalizing the matrix by multiplying the second column (on the
right) by $1-y$. The result is:
\begin{equation*}
\left(
\begin{smallmatrix}
1-z & 0 & t-1 \\
(xyx^{-1}-y)t+1-y & (xyx^{-1}-y)t^2 + (1-xyx^{-1})t+(y-1) & 0 \\
(yz-zxyx^{-1})t+y-z & (yz-zxyx^{-1}z)t^2 + (zxyx^{-1}-z)t + (1-y) & -zxyx^{-1}t^2+yt-1 \\
\end{smallmatrix}
\right)
\end{equation*}
Next we add the first column times $1-t$ and the last column times
$1-z$ to the second column. This gives us our expected free
generator:
\begin{equation*}
\left(
\begin{array}{ccc}
1-z & 0 & t-1 \\
(xyx^{-1}-y)t+1-y & 0 & 0 \\
(yz-zxyx^{-1})t+y-z & 0 & -zxyx^{-1}t^2+yt-1 \\
\end{array}
\right)
\end{equation*}
Now we subtract the first row from the third, add the second row to
the third, and then multiply the third row (on the left) by
$t^{-1}$. The result is:
\begin{equation*}
\left(
\begin{array}{ccc}
1-z & 0 & t-1 \\
(xyx^{-1}-y)t+1-y & 0 & 0 \\
x^{-1}yxz-zy+y-x^{-1}yx & 0 & -zyt+x^{-1}yx-1 \\
\end{array}
\right)
\end{equation*}
Next we add $zy$ times the first row to the third:
\begin{equation*}
\left(
\begin{array}{ccc}
1-z & 0 & t-1 \\
(xyx^{-1}-y)t+1-y & 0 & 0 \\
x^{-1}yxz+y-x^{-1}yx-zyz & 0 & x^{-1}yx-1-zy \\
\end{array}
\right)
\end{equation*}

We now have to consider two cases: whether or not $z \in
(\mathfrak{B}_4)_r^{(n)}$. From \cite{MR}, we know that $z \in
(\mathfrak{B}_4)_r^{(3)}$, but it is unclear if this holds for
$n\geq4$. If $z \in (\mathfrak{B}_4)_r^{(n+1)}$, then $z=1$ in
$\Z\Gamma_n$. In this case, our presentation matrix is:
\begin{equation*}
\left(
\begin{array}{ccc}
0 & 0 & t-1 \\
(xyx^{-1}-y)t+1-y & 0 & 0 \\
0 & 0 & x^{-1}yx-1-y \\
\end{array}
\right)
\end{equation*}
Since $x^{-1}yx-1-y$ has three terms, it cannot be equal to zero in
$\K_n[t^{\pm1}]$, and therefore is a unit. Hence we can eliminate
the last column and row. Therefore, if $z \in
(\mathfrak{B}_4)_r^{(n+1)}$,
$$H_1(\U,u_0;\K_n[t^{\pm1}])\cong\K_n[t^{\pm1}]\oplus\K_n[t^{\pm1}]/\langle
(xyx^{-1}-y)t+1-y\rangle.$$ From \cite{MR}, we know that $y\neq
xyx^{-1}$ in $\mathfrak{B}_4^\prime/\mathfrak{B}_4^{\prime\prime}$,
and therefore $xyx^{-1}-y\neq0$ in $\K_n[t^{\pm1}]$ for $n\geq1$.
Thus, if $z \in (\mathfrak{B}_4)_r^{(n+1)}$, it follows that
$\delta_n=1$. In particular, $\delta_2=1$.

Now we consider the case where $z \notin (\mathfrak{B}_4)_r^{(n)}$.
In this case, $1-z$ is invertible in $\K_n[t^{\pm1}]$. Continuing
with our calculation above, we can then multiply the first row by
$(1-z)^{-1}$ to obtain:
\begin{equation*}
\left(
\begin{array}{ccc}
1 & 0 & (1-z)^{-1}(t-1) \\
(xyx^{-1}-y)t+1-y & 0 & 0 \\
x^{-1}yxz+y-x^{-1}yx-zyz & 0 & x^{-1}yx-1-zy \\
\end{array}
\right)
\end{equation*}
Next we multiply the first row by $(y-xyx^{-1})t+y-1$ and add it to
the second. Also we multiply the first row by
$x^{-1}yx+zyz-x^{-1}yxz-y$ and add it to the third. This allows us
to eliminate the first column and row. The result is:
\begin{equation*}
\left(
\begin{array}{cc}
0 &  (y-xyx^{-1})t(1-z)^{-1}(t-1)+(y-1)(1-z)^{-1}(t-1)\\
0 & x^{-1}yx-1-zy \\
\end{array}
\right)
\end{equation*}
Since $x^{-1}yx-1-zy$ has three terms, it cannot be equal to zero,
and therefore is a unit in $\K_n[t^{\pm1}]$. Hence,
$H_1(\U,u_0;\K_n[t^{\pm1}])\cong\K_n[t^{\pm1}]$. Thus if $z \notin
(\mathfrak{B}_4)_r^{(n+1)}$, $\delta_n=0$.

To summarize, for curves whose complement has the fundamental group
$\mathfrak{B}_4$, we have shown that $\delta_0=2$ and $\delta_1=1$.
Furthermore, $\delta_n=1$ as long as $z \in
(\mathfrak{B}_4)_r^{(n+1)}$. If this is not the case, then
$\delta_n=0$. So if it can be shown that $z \notin
(\mathfrak{B}_4)_r^{(\omega)}$, then there is an integer $m\geq2$
such that $\delta_n=0$ for all $n\geq m$.

The same kind of calculations can be carried out for the other Artin
groups of spherical-type. We summarize these results without
providing the explicit calculations. For the Artin group of type
$A_2$, $\delta_0=2$ and $\delta_n=1$ for $n\geq1$. For the Artin
group of type $B_2$, $\delta_n=2$ for $n\geq0$. For the Artin group
of type $B_3$, $\delta_0=4$ and $\delta_n=3$ for $n\geq1$. For the
Artin groups of types $B_4$ and $F_4$, $\delta_0=1$ and $\delta_n=0$
for $n\geq1$. For the Artin group of type $D_4$, $\delta_0=2$,
$\delta_1=1$, and $\delta_n=0$ for $n\geq2$. For the Artin group of
type $I_2(m)$, where $m$ is odd, $\delta_0=m-1$ and $\delta_n=m-2$
for $n\geq1$. For the Artin group of type $I_2(m)$, where $m$ is
even, $\delta_n=m-2$ for $n\geq0$.

\subsection{Obstructions on fundamental groups of plane curve complements}
It follows from Theorem \ref{f}, that if $C$ is a plane curve in
general position at infinity, then $\delta_n<\infty$ for all
$n\geq0$. This is not true for a free group with at least three
generators. Therefore, such a group cannot be the fundamental group
of a plane curve complement \emph{in general position at infinity}.
Also Harvey \cite{H} has shown that $\delta_0=\infty$ for the
fundamental group of a boundary link complement. (Recall that a
boundary link is a link whose components bound mutually disjoint
Seifert surfaces.) An example of such a group (which is Ex. 8.3 of
\cite{H}) is:
\begin{eqnarray*}&&\langle a, b, c, d, e, f, g, h, i, j, k, l \;|\;
bg^{-1}ic^{-1}i^{-1}g, cj^{-1}la^{-1}l^{-1}j,
fe^{-1}hg^{-1}h^{-1}e,\\&&\hspace{.5in} ih^{-1}kj^{-1}k^{-1}h,
lk^{-1}ed^{-1}e^{-1}k, da^{-1}e^{-1}a, ebf^{-1}b^{-1},
gb^{-1}h^{-1}b,\\&&\hspace{.5in} hci^{-1}c^{-1}, jc^{-1}k^{-1}c,
kal^{-1}a^{-1}\rangle\end{eqnarray*} Therefore such a group cannot
be the fundamental group of a plane curve complement in general
position at infinity.  These facts can be combined with the sequence
(\ref{les}) in order to obtain classes of groups that cannot be the
fundamental group of a projective plane curve complement.

\providecommand{\bysame}{\leavevmode\hbox
to3em{\hrulefill}\thinspace}

\end{document}